\documentclass[11pt]{amsart}
\usepackage[english]{babel}
\usepackage{amssymb}
\usepackage{enumerate}
\usepackage{amsthm}
\usepackage{amsmath} 
\usepackage[all,cmtip]{xy}

\usepackage{tikz}
\usetikzlibrary{decorations.pathreplacing}

\usepackage[margin=1in]{geometry}
\usepackage{setspace}

\usepackage{hyperref}

\newtheorem{theorem}{Theorem}
\newtheorem{corollary}[theorem]{Corollary}

\newtheorem{lemma}[theorem]{Lemma}

\newtheorem{remark}[theorem]{Remark}

\numberwithin{theorem}{section}

\DeclareMathOperator{\area}{area} 
\DeclareMathOperator{\vol}{vol}
\newcommand{\R}{\mathbb{R}}
\DeclareMathOperator{\dist}{dist}

\title{Area-minimizing surfaces in asymptotically flat three-manifolds}

\author{Otis Chodosh}
\address{Department of Pure Mathematics and Mathematical Statistics, University of Cambridge, Wilberforce Road, Cambridge CB3 0WB, United Kingdom}
\email{oc249@cam.ac.uk}

\author{Michael Eichmair}
\address{Faculty of Mathematics, University of Vienna, Oskar-Morgenstern-Platz 1, 1090 Vienna, Austria}
\email {michael.eichmair@univie.ac.at}

\begin{document}

\maketitle

\begin{abstract} 
We show that an asymptotically flat Riemannian three-manifold with non-negative scalar curvature is isometric to flat $\R^3$ if it admits an unbounded area-minimizing surface. This answers a question of R. Schoen.
\end{abstract}

\section {Introduction}

A Riemannian three-manifold $(M, g)$ is \emph{asymptotically flat} if it is connected and there is a coordinate chart \emph{at infinity} 
\[
M \setminus U \cong \{x \in \R^3 : |x| > 1\}
\]
where $U \subset M$ is bounded and open such that 
\begin{eqnarray} \label{eqn:AF}
g_{ij} = \delta_{ij} + O(|x|^{-p})
\end{eqnarray}
and 
\[
\partial_\ell g_{ij} = O (|x|^{-p - 1}) \quad \text{ and } \quad  \partial_k \partial_\ell g_{ij} = O(|x|^{-p - 2})
\]
for some $p > 1/2$. We also require that the scalar curvature of $(M, g)$ is integrable and that the boundary of $M$ is a minimal surface in $(M,g)$. Given $\rho > 1$, we denote by $S_\rho$ the surface in $M$ corresponding to the coordinate sphere $\{x \in \R^3 : |x| = \rho\}$, and by $B_\rho$ the bounded open subset of $M$ that is enclosed by $S_\rho$. The \emph{ADM-mass} (after R. Arnowitt, S. Deser, and C. Misner \cite{ADM:1961}) of such $(M, g)$ is defined by
\[
m_{ADM} = \lim_{\rho \to \infty} \frac{1}{16 \pi}  \int_{\{|x| = \rho\}} \sum_{i, j = 1}^3 \left( \partial_i g_{ij} - \partial_j g_{ii} \right) \frac{x^j}{\rho}
\]
where the integral is with respect to the Euclidean metric. It was shown by R. Bartnik \cite{Bartnik:1986} that this definition is independent of the particular choice of the chart at infinity. 

It would be hard to overstate the significance of the \emph{positive mass theorem} as well as that of its proof using minimal surface theory  due to R. Schoen and S.-T. Yau.

\begin {theorem} [Positive mass theorem \cite{PMT1}] Let $(M, g)$ be an asymptotically flat Riemannian three-manifold with non-negative scalar curvature. Then $m_{ADM} \geq 0$. Equality holds if and only if $(M, g)$ is isometric to $\R^3$ with the flat metric. 
\end {theorem}

The proof that $m_{ADM} \geq 0$ in \cite{PMT1} is based on the insight that asymptotically flat Riemannian three-manifolds with positive scalar curvature do not admit planar two-sided stable minimal surfaces with quadratic area growth. A review of the full argument that is tailored to our discussion here can be found in the introduction of \cite{effectivePMT}.

The proof of the positive mass theorem suggests the following conjecture of R. Schoen: \emph{An asymptotically flat Riemannian manifold with non-negative scalar curvature that contains an unbounded area-minimizing surface is isometric to Euclidean space.} Cf. \cite[p. 48]{Schoen:talk}. The goal of this paper is to prove this conjecture:

\begin {theorem} \label{thm:main} The only asymptotically flat Riemannian three-manifold with non-negative scalar curvature that admits a non-compact area-minimizing boundary is flat $\R^3$. 
\end {theorem}

We make precise our meaning of \emph{area-minimizing boundaries} in Appendix \ref{sec:am}. The following surprising result of A. Carlotto and R. Schoen shows that the condition that $\Sigma \subset M$ be area-minimizing is sharp.

\begin {theorem} [A. Carlotto and R. Schoen \cite{Carlotto-Schoen}] \label{thm:CS} There exists an asymptotically flat Riemannian metric $g = g_{ij} \, dx^i \otimes dx^j$ with non-negative scalar curvature and positive mass on $\R^3$ such that $g_{ij} = \delta_{ij}$ on $\R^2 \times (0, \infty)$.
\end {theorem} 

In particular, the coordinate planes $\R^2 \times \{z\}$ with $z > 0$ in Theorem \ref{thm:CS}  are two-sided stable minimal surfaces. These planes are \emph{not} area-minimizing by Theorem \ref{thm:main}. 

On the other hand, when $(M, g)$ is \emph{asymptotic to Schwarzschild} with mass $m > 0$ in that instead of \eqref{eqn:AF} we require that  
\[
g_{ij} = \left(1 + \frac{m}{2 |x|} \right)^4 \delta_{ij} + o (|x|^{-1})
\]
with corresponding estimates for the derivatives, then Theorem \ref{thm:main} follows from the work of A. Carlotto \cite{Carlotto:2014}. We also point out the following rigidity result for two-sided stable minimal immersions. 

\begin {theorem} [A. Carlotto, O. Chodosh, and M. Eichmair \cite{effectivePMT}] \label{thm:effectivePMT} Let $(M, g)$ be a Riemannian three-manifold with non-negative scalar curvature that is asymptotic to Schwarzschild with $m > 0$. There is no unbounded complete two-sided stable minimal immersion $\varphi : \Sigma \to M$ that does not cross itself.
\end {theorem} 

We refer the reader to the introduction of our paper \cite{effectivePMT} with A. Carlotto for further remarks and a survey of related results. Our proof of Theorem \ref{thm:main} is of a very different flavor. 

Theorem \ref{thm:main} has an important consequence for the study of large isoperimetric regions. In \cite{effectivePMT}, we observe that the ideas of Y. Shi \cite{Shi:isoIMCF} imply the following existence result: given an asymptotically flat Riemannian three-manifold $(M, g)$ with non-negative scalar curvature and $V > 0$, among all smooth compact regions $\Omega \subset M$ with $\partial M \subset \Omega$ and $\vol(\Omega) = V$, there is one whose boundary has least area. Let $\Omega_V$ denote one such \emph{isoperimetric region} of volume $V$. The boundary of $\Omega_V$ is a closed stable constant mean curvature surface. J. Metzger and the second named author have shown in \cite{isostructure} that $\Omega_V$ is unique when $(M, g)$ is asymptotic to Schwarzschild with mass $m > 0$ and $V>0$ is sufficiently large. Moreover, its ``outer boundary" is close to $S_\rho$ where $\rho > 1$ is such that $\text {vol}(B_\rho) = V$. The isoperimetric regions of the \emph{exact} Schwarzschild geometry 
\[
g_{ij} = \left( 1 + \frac{m}{2 |x|}\right)^4 \delta_{ij} \quad \text{ on } \quad \{x \in \R^3 : |x| \geq m/2\}
\]
where $m > 0$ had been characterized in the fundamental work of H. Bray \cite{Bray:1997}. We also mention the seminal uniqueness results of G. Huisken and S.-T. Yau \cite{Huisken-Yau:1996} and of J. Qing and G. Tian \cite{Qing-Tian:2007} for large two-sided stable constant mean curvature spheres in Riemannian three-manifolds that are asymptotic to Schwarzschild with mass $m > 0$, as well as extensions of these results obtained in \cite{stableCMC, effectivePMT, offcenter}. 

It is natural to wonder about the behavior of large isoperimetric regions when $(M, g)$ has general asymptotics. For simplicity of exposition, let us assume that the boundary of $M$ is empty. Let $\Sigma_i = \partial \Omega_{V_i}$ be a sequence of isoperimetric surfaces with $0 < V_i \to \infty$. It has been shown in \cite{isostructure} that these surfaces either diverge to infinity as $i \to \infty$, or that   a subsequence of these surfaces converges geometrically to a non-compact area-minimizing boundary $\Sigma \subset M$. In view of Theorem \ref{thm:main}, the latter is impossible unless $(M, g)$ is flat $\R^3$. In conclusion, we arrive at the dichotomy that large isoperimetric regions in $(M, g)$ are either drawn far into the asymptotically flat end, or they contain the center of the manifold. When the scalar curvature of $(M, g)$ is everywhere positive, this was observed in Corollary 6.2 of \cite{isostructure}.

\begin {corollary} Let $(M, g)$ be an asymptotically flat Riemannian three-manifold with non-negative scalar curvature and positive mass. We also assume that the components of the boundary of $M$ are the only closed minimal surfaces in $(M, g)$. Let $U \subset M$ be a bounded open subset that contains the boundary of $M$. There is $V_0 > 0$ so that for every isoperimetric region $\Omega \subset M$ of volume $V \geq V_0$, either $U \subset \Omega_V$ or $U \cap \Omega_V$ is a thin  smooth region that is bounded by the components of $\partial M$ and nearby two-sided stable constant mean curvature surfaces.
\end {corollary} 
Note that the conclusion of the corollary is wrong for flat $\R^3$. The role of Theorem \ref{thm:main} in the proof of this corollary is similar to that of Theorem \ref{thm:effectivePMT} in the analysis \cite{effectivePMT} of large stable constant mean curvature surfaces in three-manifolds that are asymptotic to Schwarzschild with mass $m>0$.

Finally, we mention several other results that are based on scalar curvature and the existence of two-sided stable minimal surfaces of a certain topological type. 

\begin {theorem} [R. Schoen and S.-T. Yau \cite{Schoen-Yau:1978}] A Riemannian three-manifold with positive scalar curvature does not admit a two-sided stable minimal immersion of positive genus. 
\end {theorem}

\begin{theorem}[M. Cai and G. Galloway \cite{CaiGalloway}]
Let $(M,g)$ be a closed Riemannian three-manifold with non-negative scalar curvature. If $(M,g)$ contains an embedded area-minimizing two-torus, then $(M,g)$ is flat.
\end{theorem}

\begin {theorem} [H. Bray, S. Brendle, M. Eichmair, and A. Neves \cite{Bray-Brendle-Eichmair-Neves:2010}] Let $(M, g)$ be a closed Riemannian three-manifold that contains an embedded projective plane. We have that
\[
 \inf \{ \area(\Sigma) : \Sigma \subset M \text{ is an embedded projective plane}\}  \cdot  \inf_M R  \leq 12 \pi
\]
where $R$ denotes the scalar curvature of $(M, g)$. Equality holds if and only if $(M, g)$ is isometric to real projective space $\mathbb{R P}^3$ with the round metric, up to scaling. 
\end {theorem}

There have recently been some important developments in this area. We refer the reader to \cite{Bray-Brendle-Neves:2010,MarquesNeves:min-max-rigidity-3mflds,MaximoNunes:hawking-rigidity,Nunes:hyperbolic-rigidity,Ambrozio:free-bdry-rigidity,MicallefMoraru} and the references therein. \\

{\bf Acknowledgments: } We are very grateful to H. Bray, S. Brendle, G. Galloway, J. Metzger, and R. Schoen for their support. It is a pleasure to sincerely congratulate R. Schoen on the occasion of his 65th birthday. 


\section {Tools}

The following result is a technical variation of a crucial insight in the proof of the positive mass theorem due to R. Schoen and S.-T. Yau  \cite{PMT1}. 

\begin {lemma} [Section 6 in \cite{isostructure}] \label{lem:simplepmt} Let $(M, g)$ be an asymptotically flat Riemannian three-manifold. Assume that $\Sigma \subset M$ is the unbounded component of an area-minimizing boundary in $(M, g)$, and that the scalar curvature of $(M, g)$ is non-negative along $\Sigma$. Then $\Sigma \subset M$ is totally geodesic and the   scalar curvature of $(M, g)$ vanishes along this surface. Moreover, for all $\rho > 1$ sufficiently large, $\Sigma$ intersects $S_\rho$ transversely in a nearly equatorial circle.
\end {lemma}

We also need the following rigidity result for asymptotically flat slabs with minimal boundary from our paper \cite{effectivePMT}. The proof adapts to the setting of non-negative \emph{scalar} curvature a method that was developed by M. Anderson and L. Rodr\'iguez \cite{Anderson-Rodriguez:1989} and crucially refined by G. Liu \cite{Liu:2013} to characterize complete manifolds with non-negative \emph{Ricci} curvature. 

\begin {lemma} [Theorem 1.5 in \cite{effectivePMT}] \label{lem:slab} 
Let $(M, g)$ be an asymptotically flat Riemannian three-manifold with non-negative scalar curvature. Any two disjoint unbounded connected complete properly embedded minimal surfaces in $(M,g)$ bound a region that is isometric to a Euclidean slab $\R^2 \times [a,b]$.\footnote{This result is stated in \cite{effectivePMT} with the additional assumption that the components of the boundary of $M$ are the only closed minimal surfaces in $(M,g)$. By choosing the curve $\gamma$ in the proof of Theorem 1.5 in \cite{effectivePMT} to lie sufficiently far out, we ensure that every complete minimal surface that intersects $\gamma$ is unbounded. A similar argument is made in the choice of the point $p$ in Section \ref{sec:proof}. The rest of the proof in \cite{effectivePMT} goes through verbatim.}
\end {lemma}


\section {Proof of Theorem \ref{thm:main}} \label{sec:proof}

We first deal with the case where the boundary of $M$ is empty. 

Let $r_0 > 0$ be as in Appendix \ref{sec:cf}. Let $\rho_0 > 1$ be such that $S_\rho$ is convex for all $\rho \geq \rho_0$. Every closed minimal surface of $(M, g)$ is contained in $B_{\rho_0}$. 

Let $\Sigma \subset M$ be a connected unbounded properly embedded and separating surface that is area-minimizing with respect to $g$. Fix a component $M_+$ of the complement of $\Sigma$ in $M$ and choose a point $p \in M_+$ to the following specifications: 
\begin {itemize}
\item $r = \dist_g(x, p)/2 < r_0$;
\item $p \notin B_{\rho_0 + 4 r_0}$; 
\item $\Sigma$ intersects $\{x \in M : \dist_g(x, p) < 4 r\}$ in a single component, and this component is transverse to $\nabla_g \dist_g(\, \cdot \, , p)$.
\end {itemize}

\begin{figure}[h!] 
\begin{tikzpicture}
	\filldraw [white] (-7,-1.2) rectangle (7,3.4);

	\filldraw [opacity = .2] (-2.5,1) circle (1.8);
	\filldraw [white] (-2.5,1) circle (.6);
	
	\draw (-2.5,1) circle (.6);
	\draw (-2.5,1) circle (1.8);
	
	\filldraw (-2.5,1) circle (1pt) node [shift={(-.2,-.1)}] {$p$};
	\draw [->] (-2.5,1) -- node [shift = {(.1,.15)}] {$\scriptstyle r$} +(-33:.6);
	\draw [->] (-2.5,1) -- node [shift = {(.2,-.25)}] {$\scriptstyle 2r$} +(-93:1.1);
	\draw [->] (-2.5,1) -- node [shift = {(.35,.2)}] {$\scriptstyle 3r$} +(70:1.8);

	\begin{scope}
		\clip (-7,-1.2) rectangle (7,3.2);
		\draw [dashed] (0,0) circle (6);
	\end{scope}
	
	\draw plot [smooth] coordinates {(-7,-.3) (-3,-.1) (3,-.5) (7,0)};
	\node at (3,-.9) {$\Sigma$};
	
	\draw plot [smooth] coordinates {(182:6) (-1,1.5)  (-1.7:6)};
	\node at (3,1.2) {$\Sigma_{\rho}(t)$};
	
	\node at (-6.3,-.6) {$\Gamma_{\rho}$};
	\node at (6,2) {$S_{\rho}$};
	
	\node at (1.5,2.5) {$\scriptstyle g(t) = g$};
	\node at (-3.3,1.8) {$\scriptstyle R_{g(t)} > 0$};

	\node at (-6.7,2.3) {$M_{+}$};
\end{tikzpicture}
\caption{A diagram of the perturbed metric $g(t)$ and corresponding surface $\Sigma_{\rho}(t)$ used in the proof of Theorem \ref{thm:main}.} \label{Figure 1}
\end{figure}
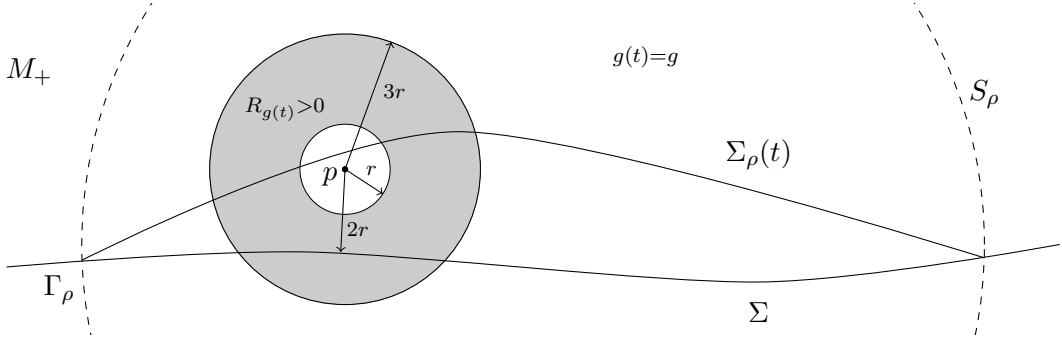

In Appendix \ref{sec:cf}, we construct a family of conformal Riemannian metrics $\{g(t)\}_{t \in [0, \epsilon)}$ on $M$ with the following properties (see also Figure \ref{Figure 1}): 
\begin {enumerate} [(i)]
\item $g(t) \to g$ smoothly as $t \to 0$;
\item $g(t) = g$ on $\{x \in M :  \dist_g (x, p) \geq 3 r \}$;
\item $g(t) \leq g$ as quadratic forms on $M$, with strict inequality on $\{x \in M : r < \dist_g (x, p) < 3 r\}$;
\item the scalar curvature of $g(t)$ is positive on $\{x \in M : r < \dist_g(x, p) < 3 r\}$;
\item the region $M_+$ is weakly mean-convex with respect to $g(t)$. 
\end {enumerate}

By taking $\epsilon >0$ smaller if necessary, we may assume that all closed minimal surfaces of $(M, g(t))$ are contained in $B_{\rho_0}$.

According to Lemma \ref{lem:simplepmt}, for all $\rho \geq \rho_0$ sufficiently large, the intersection of $\Sigma$ with $S_\rho$ is transverse in a nearly equatorial circle. We denote this circle by $\Gamma_\rho = \Sigma \cap S_\rho$. Consider all properly embedded surfaces in $M$ that have boundary $\Gamma_\rho$ and which together with $\Sigma \cap B_\rho$ bound an open subset of $M_+ \cap B_\rho$. Using (v) and standard existence results from geometric measure theory, we see that among all these surfaces there is one  --- call it $\Sigma_\rho (t)$ ---  that has least area with respect to $g(t)$. This surface is disjoint from $M_+ \cap S_\rho$ by convexity. It has one component with boundary $\Gamma_\rho$. Its other components are closed minimal surfaces in $(M, g(t))$. In particular, they are disjoint from $\{x \in M : \dist_g(x, p) < 3 r\}$. Importantly, $\Sigma_\rho(t)$ \emph{does intersect} $\{x \in M : \dist_g (x, p) < 3 r\}$ since otherwise,
\begin{eqnarray} \label{eqn:otherwise}
\area_g (\Sigma_\rho(t)) = \area_{g(t)} (\Sigma_\rho(t)) \leq \area_{g(t)} (\Sigma \cap B_\rho) < \area_g (\Sigma \cap B_\rho).
\end{eqnarray}
The strict inequality holds on account of (iii) and because $\Sigma$ intersects $\{x \in M : \dist_g(x, p) < 3 r\}$. Observe that \eqref{eqn:otherwise} violates the area-minimizing property of $\Sigma$ with respect to $g$.

Using standard convergence results from geometric measure theory, we now find a connected unbounded properly embedded separating surface $\Sigma(t) \subset M$ as a subsequential geometric limit of $\Sigma_\rho(t)$ as $\rho \to \infty$. By construction, $\Sigma(t)$ is contained in $M_+ \cup \Sigma$ where it is area-minimizing with respect to $g(t)$. Moreover, $\Sigma (t)$ intersects $\{x \in M : \dist_g (x, p) \leq 3 r\}$. If $\Sigma (t)$ intersects $\{x \in M : \dist_g(x, p) < 3 r\}$, then it also intersects $\{x \in M : \dist_g(x, p) \leq r\}$ because of (iv) and Lemma \ref{lem:simplepmt}. Passing to a subsequential geometric limit as $t \to 0$, we obtain a connected unbounded properly embedded separating surface $\Sigma_+ \subset M$ that is contained in $M_+ \cup \Sigma$ where it is area-minimizing with respect to  $g(t)$. Using now the area-minimizing property of $\Sigma$, we see that $\Sigma_+$ is in fact area-minimizing in all of $M$. Note that $\Sigma$ intersects $\{x \in M : \dist_g(x, p) < 3 r\}$ while it is disjoint from $\{x \in M : \dist_g(x, p) \leq r\}$. It follows from the maximum principle that $\Sigma$ and $\Sigma_+$ are disjoint. The region in $(M, g)$ that is bounded by $\Sigma$ and $\Sigma_+$ is isometric to a Euclidean slab by Lemma \ref{lem:slab}. 

Consider now the set of all $h > 0$ so that $(M, g)$ contains a Euclidean slab of the form $\R^2 \times [0, h]$ with $\R^2 \times \{0\}$ corresponding to $\Sigma$ and $\R^2 \times \{h\}$ corresponding to another area-minimizing surface in $(M, g)$. We have seen that this set is non-empty and open. It is clearly closed. It follows that $(M, g)$ is isometric to a Euclidean half-space to one side (namely $M_+$) of $\Sigma$. Since the choice of $M_+$ as one of the two components of the complement of $\Sigma$ in $M$ was arbitrary, we conclude that $(M, g)$ is indeed isometric to flat $\R^3$.

We now turn to the general case where $M$ has boundary. Consider $\Omega \in \mathcal{F}$ with non-compact area-minimizing boundary $\Sigma \subset M$. The unique non-compact component $\Sigma_0 \subset M$ of $\Sigma$ is a separating surface. Let $M_-$ and $M_+$ denote the two components of its complement in $M$. Note that the interior of $\Omega \cap M$ agrees with either $M_-$ (Case 1) or $M_+$ (Case 2) outside of $B_{\rho_0}$. The proof that $g$ is flat in $M_+$ proceeds exactly as above, except for the following change. In Case 1, we let $\Sigma_\rho (t)$ have least area among properly embedded surfaces with boundary $\Gamma_\rho$ that bound together with $\Sigma_0 \cap B_\rho$ in $M_+ \cap B_\rho$ and relative to $M_+ \cap \partial M$. In Case 2, we let $\Sigma_\rho(t)$ have least area among properly embedded surfaces with boundary $\Gamma_\rho$ that bound together with $M_+ \cap S_\rho$ in $M_+ \cap B_\rho$ and relative to $M_+ \cap  \partial M$. Theorem \ref{thm:main} follows upon switching the roles of $M_-$ and $M_+$.

\begin {remark} The use of the conformal change of metric in this proof is inspired by an idea of G. Liu in his classification of complete non-compact Riemannian three-manifolds with non-negative Ricci curvature \cite{Liu:2013}. This idea has partially been adapted to non-negative scalar curvature in our proof of Lemma \ref{lem:slab}. The observation \eqref{eqn:otherwise} in the proof of Theorem \ref{thm:main} is crucial: we use it to conclude that the surfaces $\Sigma_\rho(t)$ cannot run off as $\rho \to \infty$. At a related point in the work of M. Anderson and L. Rodr\'iguez \cite{Anderson-Rodriguez:1989}, the assumption of non-negative Ricci curvature is used tacitly in their delicate estimation of comparison surfaces \cite[(1.5)]{Anderson-Rodriguez:1989}. 
\end {remark} 


\appendix

\section {} \label{sec:am}

Let $(M, g)$ be an asymptotically flat Riemannian three-manifold. Extend $M$ inwards across each of its minimal boundary components by thin open collar neighborhoods to a new manifold $\hat M$. Denote the union of these finitely many collar neighborhoods by $C$. We consider the collection $\mathcal{F}$ of all proper three-dimensional submanifolds with boundary $\Omega \subset \hat M$ with $C \subset \Omega$. A surface $\Sigma \subset M$ \emph{bounds} in $M$ \emph{relative to} $\partial M$ if it arises as the boundary of such a smooth region. Note that $\partial M$ bounds in this sense. 

We say that the boundary of $\Omega \in \mathcal{F}$ is \emph{area-minimizing} if for all $\rho > 1$ and $\tilde \Omega \in \mathcal{F}$ with $ \Omega \setminus B_\rho = \tilde \Omega \setminus B_\rho$ we have that
\[
\area (B_{2\rho} \cap \partial \Omega) \leq \area(B_{2\rho} \cap \partial  \tilde \Omega).
\]
The components of the boundary $\Sigma \subset M$ of such $\Omega \in \mathcal{F}$ are two-sided stable minimal surfaces in $(M, g)$. It follows from the arguments in Section 6 of \cite{isostructure} that $\Sigma$ has at most one unbounded component $\Sigma_0$. More precisely, if we consider the homothetic blow-downs of $\Omega$ in the chart at infinity 
\[
M \setminus U \cong \{x \in \R^3 : |x| > 1\}
\]
by a sequence $\lambda_i \to \infty$, then we can pass to a geometric subsequential limit in $\R^3 \setminus \{0\}$. This limit is either a half-space through the origin, or empty in the case where $\Omega$ is bounded. 

We mention that when $M$ is orientable and when there are no closed minimal surfaces in $(M, g)$ other than the components of the boundary of $M$, then $M$ is diffeomorphic to the complement in $\R^3$ of the union of finitely many open balls whose closures are disjoint; see Section 4 in \cite{Huisken-Ilmanen:2001} and the references therein.


\section {} \label{sec:cf}

Let $f \in C^\infty(\R)$ be a non-positive function with support in the interval $[0, 3]$ such that 
\[
f(s) = - \exp (18 / (s-3))
\]
when $s \in (1, 3)$. This definition is made so that 
\[
0 <  f'(s) \quad \text{ and } \quad s f'' (s) + 3 f' (s) < 0 
\]
for all $s \in (1, 3)$. 

Let $(M, g)$ be a homogeneously regular Riemannian three-manifold. Choose $0 < r_0 < \text{inj} (M, g)/4$ so that 
\[
\Delta \dist_g(\, \cdot\, , p)^2 \leq 8 \quad \text{ on } \quad \{ x \in M : \dist_g(x, p) \leq 3 r_0\}
\] 
for all $p \in M$. Fix $p \in M$ and $0 < r \leq r_0$. Consider the function $v : M \to \R$ given by
\[
x \mapsto f ( \dist_g(x, p)/r).
\]
Note that $v$ is smooth, non-positive, and supported in $\{x \in M : \dist_g(x, p) \leq 3 r\}$. Moreover, 
\[
v< 0, \qquad g (\nabla_g v , \nabla_g \dist_g( \, \cdot\, , p))>0,  \qquad \Delta_g v < 0
\]
on $\{ x \in M : r < \dist_g(x, p) < 3 r \}$. 
For $\epsilon > 0$ sufficiently small, a smooth family of conformal metrics $\{g(t)\}_{t \in [0, \epsilon)}$ with the properties needed in the proof of Theorem \ref{thm:main} is given by
\[
g(t) = (1 + t v)^4 g. 
\]

\bibliography{bib} 

\providecommand{\bysame}{\leavevmode\hbox to3em{\hrulefill}\thinspace}
\providecommand{\MR}{\relax\ifhmode\unskip\space\fi MR }
\providecommand{\MRhref}[2]{%
  \href{http://www.ams.org/mathscinet-getitem?mr=#1}{#2}
}
\providecommand{\href}[2]{#2}
\begin{thebibliography}{10}

\bibitem{Ambrozio:free-bdry-rigidity}
Lucas Ambrozio, \emph{Rigidity of area-minimizing free boundary surfaces in
  mean convex three-manifolds}, J. Geom. Anal. \textbf{25} (2015), no.~2,
  1001--1017. \MR{3319958}

\bibitem{Anderson-Rodriguez:1989}
Michael Anderson and Lucio Rodr{\'{\i}}guez, \emph{Minimal surfaces and
  {$3$}-manifolds of nonnegative {R}icci curvature}, Math. Ann. \textbf{284}
  (1989), no.~3, 461--475. \MR{1001714}

\bibitem{ADM:1961}
Richard Arnowitt, Stanley Deser, and Charles Misner, \emph{Coordinate
  invariance and energy expressions in general relativity.}, Phys. Rev. (2)
  \textbf{122} (1961), 997--1006. \MR{0127946}

\bibitem{Bartnik:1986}
Robert Bartnik, \emph{The mass of an asymptotically flat manifold}, Comm. Pure
  Appl. Math. \textbf{39} (1986), no.~5, 661--693. \MR{849427}

\bibitem{Bray-Brendle-Eichmair-Neves:2010}
H.~Bray, S.~Brendle, M.~Eichmair, and A.~Neves, \emph{Area-minimizing
  projective planes in 3-manifolds}, Comm. Pure Appl. Math. \textbf{63} (2010),
  no.~9, 1237--1247. \MR{2675487}

\bibitem{Bray:1997}
Hubert Bray, \emph{The {P}enrose inequality in general relativity and volume
  comparison theorems involving scalar curvature}, ProQuest LLC, Ann Arbor, MI,
  1997, Thesis (Ph.D.)--Stanford University. \MR{2696584}

\bibitem{Bray-Brendle-Neves:2010}
Hubert Bray, Simon Brendle, and Andre Neves, \emph{Rigidity of area-minimizing
  two-spheres in three-manifolds}, Comm. Anal. Geom. \textbf{18} (2010), no.~4,
  821--830. \MR{2765731}

\bibitem{offcenter}
Simon Brendle and Michael Eichmair, \emph{Large outlying stable constant mean
  curvature spheres in initial data sets}, Invent. Math. \textbf{197} (2014),
  no.~3, 663--682. \MR{3251832}

\bibitem{CaiGalloway}
Mingliang Cai and Gregory Galloway, \emph{Rigidity of area minimizing tori in
  3-manifolds of nonnegative scalar curvature}, Comm. Anal. Geom. \textbf{8}
  (2000), no.~3, 565--573. \MR{1775139}

\bibitem{Carlotto:2014}
Alessandro Carlotto, \emph{Rigidity of stable minimal hypersurfaces in
  asymptotically flat spaces}, preprint,
  \url{http://arxiv.org/pdf/1403.6459.pdf} (2014).

\bibitem{effectivePMT}
Alessandro Carlotto, Otis Chodosh, and Michael Eichmair, \emph{Effective
  versions of the positive mass theorem}, preprint,
  \url{http://arxiv.org/abs/1503.05910} (2015).

\bibitem{Carlotto-Schoen}
Alessandro Carlotto and Richard Schoen, \emph{Localized solutions of the
  {E}instein constraint equations}, preprint,
  \url{http://arxiv.org/abs/1407.4766} (2014).

\bibitem{stableCMC}
Michael Eichmair and Jan Metzger, \emph{On large volume preserving stable {CMC}
  surfaces in initial data sets}, J. Differential Geom. \textbf{91} (2012),
  no.~1, 81--102. \MR{2944962}

\bibitem{isostructure}
\bysame, \emph{Large isoperimetric surfaces in initial data sets}, J.
  Differential Geom. \textbf{94} (2013), no.~1, 159--186. \MR{3031863}

\bibitem{Huisken-Ilmanen:2001}
Gerhard Huisken and Tom Ilmanen, \emph{The inverse mean curvature flow and the
  {R}iemannian {P}enrose inequality}, J. Differential Geom. \textbf{59} (2001),
  no.~3, 353--437. \MR{1916951}

\bibitem{Huisken-Yau:1996}
Gerhard Huisken and Shing-Tung Yau, \emph{Definition of center of mass for
  isolated physical systems and unique foliations by stable spheres with
  constant mean curvature}, Invent. Math. \textbf{124} (1996), no.~1-3,
  281--311. \MR{1369419}

\bibitem{Liu:2013}
Gang Liu, \emph{3-manifolds with nonnegative {R}icci curvature}, Invent. Math.
  \textbf{193} (2013), no.~2, 367--375. \MR{3090181}

\bibitem{MarquesNeves:min-max-rigidity-3mflds}
Fernando Marques and Andr{\'e} Neves, \emph{Rigidity of min-max minimal spheres
  in three-manifolds}, Duke Math. J. \textbf{161} (2012), no.~14, 2725--2752.
  \MR{2993139}

\bibitem{MaximoNunes:hawking-rigidity}
Davi M{\'a}ximo and Ivaldo Nunes, \emph{Hawking mass and local rigidity of
  minimal two-spheres in three-manifolds}, Comm. Anal. Geom. \textbf{21}
  (2013), no.~2, 409--432. \MR{3043752}

\bibitem{MicallefMoraru}
Mario Micallef and Vlad Moraru, \emph{Splitting of 3-manifolds and rigidity of
  area-minimising surfaces}, Proc. Amer. Math. Soc. \textbf{143} (2015), no.~7,
  2865--2872. \MR{3336611}

\bibitem{Nunes:hyperbolic-rigidity}
Ivaldo Nunes, \emph{Rigidity of area-minimizing hyperbolic surfaces in
  three-manifolds}, J. Geom. Anal. \textbf{23} (2013), no.~3, 1290--1302.
  \MR{3078354}

\bibitem{Qing-Tian:2007}
Jie Qing and Gang Tian, \emph{On the uniqueness of the foliation of spheres of
  constant mean curvature in asymptotically flat 3-manifolds}, J. Amer. Math.
  Soc. \textbf{20} (2007), no.~4, 1091--1110. \MR{2328717}

\bibitem{Schoen:talk}
Richard Schoen, \emph{Geometric questions in general relativity}, Chen-Jung Hsu
  Lecture 2, Academia Sinica, ROC, slides available at
  \url{http://w3.math.sinica.edu.tw/HSU/page/2013/Dec.3.pdf}, December 3 2013.

\bibitem{Schoen-Yau:1978}
Richard Schoen and Shing-Tung Yau, \emph{Incompressible minimal surfaces,
  three-dimensional manifolds with nonnegative scalar curvature, and the
  positive mass conjecture in general relativity}, Proc. Nat. Acad. Sci. U.S.A.
  \textbf{75} (1978), no.~6, 2567. \MR{496776}

\bibitem{PMT1}
\bysame, \emph{On the proof of the positive mass conjecture in general
  relativity}, Comm. Math. Phys. \textbf{65} (1979), no.~1, 45--76. \MR{526976}

\bibitem{Shi:isoIMCF}
Yuguang Shi, \emph{Isoperimetric inequality on asymptotically flat manifolds
  with nonnegative scalar curvature}, preprint,
  \url{http://arxiv.org/abs/1503.02350} (2015).

\end{thebibliography}
\bibliographystyle{amsplain}
\end{document}